\documentclass[twoside,11pt]{article}

\usepackage{a4,amsthm,latexsym,amsmath,amssymb,amsfonts,euscript}
\usepackage{epsfig}
\newtheorem{theorem}{{\bf{\small T}{\scriptsize HEOREM}}}[]

\newtheorem{proposition}[]{{\bf{\small P}{\scriptsize ROPOSITION}}}
\newtheorem{lemma}[]{{\bf{\small L}{\scriptsize EMMA}}}
\newtheorem{remark}[]{{\bf{\small R}{\scriptsize EMARK}}}

\renewenvironment{proof}[1]
{\noindent{{\bf{\small{ P}{\scriptsize ROOF}}}.}\hspace{0.1cm} #1} {$\;\qed$\newline}

\def\R{\mathbb R}
\def\N{\mathbb N}
\def\Z{\mathbb Z}
\newcommand{\Zd}{\mathbb Z^d}

\def\1{{\mathchoice {\rm 1\mskip-4mu l} {\rm 1\mskip-4mu l}
{\rm 1\mskip-4.5mu l} {\rm 1\mskip-5mu l}}}

\def\pee{{\mathbb P}}
\def\qee{{\mathbb Q}}

\def\E{\mathbb E}

\def\si{\sigma}
\def\la{\Lambda}

\def\D{\mathcal D}

\def\disagree{{\nleftrightarrow}}

\newcommand{\Oun}{\mathcal{O}(1)} 

\def\ldeux{\scriptscriptstyle{\ell^{2}(\N)}}
\def\ldeuxzd{\scriptscriptstyle{\ell^{2}(\Z^d)}}
\def\lun{\scriptscriptstyle{\ell^{1}(\N)}}
\def\Lq{\scriptscriptstyle{L^{q}(\pee)}}

\date{}

\begin{document}

\begin{center}
{\bf{\Large Concentration inequalities for random fields\\ via coupling}}
\end{center}

\vskip .5truecm

\centerline{{\bf J.-R. Chazottes}}
\centerline{Centre de Physique Th{\'e}orique, CNRS UMR 7644}
\centerline{F-91128 Palaiseau Cedex, France}
\centerline{{\tt jeanrene@cpht.polytechnique.fr}}
\centerline{{\bf P. Collet}}
\centerline{Centre de Physique Th{\'e}orique, CNRS UMR 7644}
\centerline{F-91128 Palaiseau Cedex, France}
\centerline{{\tt collet@cpht.polytechnique.fr}}
\centerline{{\bf C. K\"ulske}}
\centerline{Department of Mathematics and Computing Sciences}
\centerline{University of Groningen, Blauwborgje 3}
\centerline{9747 AC Groningen, The Netherlands}
\centerline{{\tt kuelske@math.rug.nl}}
\centerline{{\bf F. Redig}}
\centerline{Mathematisch Instituut Universiteit Leiden}
\centerline{Niels Bohrweg 1, 2333 CA Leiden, The Netherlands}
\centerline{{\tt redig@math.leidenuniv.nl}}

\vskip 2truecm

\begin{abstract}
We present a new and simple approach to concentration inequalities
for functions around their expectation with respect to non-product
measures, i.e., for dependent random variables.
Our method is based on coupling ideas and does not use information
inequalities.
When one has a uniform control on the coupling, this leads to
exponential concentration inequalities. 
When such a uniform control is no more possible, this leads to
polynomial or stretched-exponential concentration inequalities.   
Our abstract results apply to Gibbs random fields, in 
particular to the low-temperature Ising model which is a
concrete example of non-uniformity of the coupling. 

\bigskip 


\end{abstract}


\newpage


\section{Introduction}

By now, concentration inequalities for product measures have become
a standard and powerful tool in many areas of probability and statistics, such as
density estimation \cite{dl}, geometric probability \cite{yu}, etc.
A recent monograph about this area is \cite{ledoux} where the reader can find much more
information and relevant references.
Exponential concentration inequalities for functions of dependent,
strongly mixing random variables were
obtained for instance in \cite{Kul,marton1,marton2,marton3,samson,rio}. 
In the context of dynamical systems Collet et al. \cite{cms} obtained
an exponential concentration inequality for separately Lipschitz
functions using spectral analysis of the
transfer operator. In \cite{Kul}, C. K\"ulske obtained an exponential
concentration inequality for functions of Gibbs random fields in the
Dobrushin uniqueness regime. Therein the main input is Theorem 8.20 in
\cite{Geo} which allows to estimate uniformly the terms appearing 
in the martingale difference decomposition in terms of the Dobrushin matrix.
In \cite{marton2}, K. Marton obtained exponential concentration results
for a class of Gibbs random fields under a strong mixing condition
lying between Dobrushin-Shlosman condition and its weakening in the
sense of E. Olivieri, P. Picco and F. Martinelli.

Besides exponential concentration inequalities, polynomial concentration
inequalities easily follow from upper bounds on moments. In the 
context of product measures, bounds on the variance are well-known
\cite{devroye,dl}. In the context of dynamical systems, a bound
on the variance is obtained in \cite{ccs}. 

The approach followed in \cite{marton1,marton2,marton3,samson} uses
coupling ideas and information inequalities, such as Pinsker
inequality. Such inequalities can only lead to exponential
concentration inequalities. This can be understood easily
since it is well-known \cite{BG} that there is equivalence
between information inequalities and exponential inequalities
on the Laplace transform, the latter yielding exponential
concentration inequalities by Chebychev's inequality.

The purpose of the present paper is to derive abstract bounds allowing
to obtain not only exponential, but also polynomial and
stretched-exponential concentration inequalities. 
In particular, this means that we do not use information
inequalities. Going beyond the exponential case was motivated by the
low-temperature Ising model which can 
not satisfy an exponential concentration inequality for the magnetization.
Here we obtain abstract concentration inequalities using a coupling
approach. 
Our setting is (dependent) random variables
indexed by $\Z^d$, $d\geq 1$, and taking values in a finite alphabet.
We are interested in obtaining concentration inequalities for ``local''
functions $g$ around their expectation $\E g$ in terms of their
variations.   
The inter-dependence between random variables is
measured by a ``coupling matrix'' which tells us how ``well'' one can couple in the
far ``future'' if the ``past'' is given. If the coupling matrix can
be uniformly controlled in the realization, then an exponential
concentration inequality follows. If the coupling matrix cannot 
be controlled uniformly in the realization, then we typically obtain
bounds for moments and for Luxembourg norms of $g-\E g$. In the former
case this leads to polynomial concentration inequalities, in the latter
case this gives stretched-exponential concentration inequalities.

As a first application of our abstract inequalities,
we obtain an exponential concentration inequality for
Gibbs random fields in a ``high-temperature'' regime, complementary
to the Dobrushin uniqueness regime studied in \cite{Kul}.
A second application is the ``low-temperature'' Ising model for which
the coupling matrix cannot be uniformly controlled in the realization,
and for which the previous methods \cite{marton2,samson} do not apply. 
We obtain polynomial, even stretched-exponential, concentration
inequalities for the low-temperature Ising model. Let us
mention that our concentration inequalities yield various non-trivial
applications which will be the subject of a forthcoming paper.

The paper is organized as follows. In Section 2, we state and prove our abstract inequalities,
first in the context of random fields indexed by $\Z$, and next when
the index set is $\Z^d$, $d\geq 2$. Section 3 deals with
high-temperature Gibbs measures and the low-temperature Ising model.

\section{Main results}\label{MR}

Let $A$ be a finite set. Let $g:A^n\to\R$ be a function of $n$-variables.
An element $\si$ of the set $A^{\N}$ is an infinite sequence
drawn from $A$, i.e., $\si=(\si_1,\si_2,\ldots,\si_i,\ldots)$ where $\si_i\in A$.
With a slight abuse of notation, we also consider $g$ as a function on $A^\N$
which does not depend on $\si_{k}$, for all $k>n$.

A concentration inequality is an estimate for the probability of concentration 
of the function $g$ from its expectation, i.e., an estimate for
\begin{equation}\label{proba}
\pee\left\{ | g -\E g |\geq t\right\}
\end{equation}
for all $n\geq 1$ and all $t>0$, within a certain class of probability measures $\pee$.
For example, an {\em exponential concentration inequality} is obtained by estimating the
expectation
$$
\E\left[ e^{\lambda(g -\E g)}\right]
$$
for any $\lambda\in\R$, and using the exponential Chebychev's inequality.

However, there are natural examples where the exponential concentration inequality
does not hold (see the example of the low-temperature Ising model below).
In that case we are interested in bounding moments of the form
$$
\E\left[ (g -\E g)^{2p}\right]
$$
to control the probability \eqref{proba}.

In this section, we use a combination of the classical martingale decomposition of
$g -\E g$ and maximal coupling to perform a further telescoping
which is adequate for the dependent case. 
This will lead us to a ``coupling matrix'' depending on the realization
$\si\in A^{\N}$. This matrix quantifies how ``good'' future symbols can be coupled
if past symbols are given according to $\si$.
Typically, we have in mind applications to Gibbs random fields. In that framework,
the elements of the coupling matrix can be controlled uniformly in $\si$ in the 
``high-temperature regime''. This uniform control leads
naturally to an exponential concentration inequality.
At low temperature we can only control the coupling matrix for ``good'' configurations,
but not uniformly. Therefore an exponential concentration inequality
cannot hold (for all $g$). Instead we will obtain polynomial and
stretched-exponential concentration inequalities. This will be done by
controlling moments and Luxembourg norms of $g-\E g$.

\subsection{The coupling matrix $D^\si$}

We now present our method.
For $i=1,2,\ldots,n$, let $\mathcal{F}_i$ be the sigma-field generated by the random variables
$\si_1,\ldots,\si_i$, and $\mathcal{F}_0$ be the trivial sigma-field $\{\emptyset,\Omega\}$.
We write 
\begin{equation}\label{decomp}
g(\si_1,\ldots,\si_n) -\E g= \sum_{i=1}^{n} V_i(\si)
\end{equation}
where
$$
V_i(\si):= \E[g | \mathcal{F}_i](\si)-\E[g | \mathcal{F}_{i-1}](\si)=
$$
$$
\int \pee(d\eta_{i+1}\cdots d\eta_{n} | \si_1,\ldots,\si_i)\ g(\si_1,\ldots,\si_i, \eta_{i+1},\ldots, \eta_{n})
$$
$$
-
\int \pee(d\eta_{i}\cdots d\eta_{n} | \si_1,\ldots,\si_{i-1})\
g(\si_1,\ldots,\si_{i-1},\eta_i,\eta_{i+1},\ldots,\eta_{n})=
$$
$$
\int \pee(d\eta_{i+1}\cdots d\eta_{n} | \si_1,\ldots,\si_i)\ g(\si_1,\ldots,\si_i, \eta_{i+1},\ldots, \eta_{n})
$$
$$
-
\int \pee(d\eta_i|\si_1,\ldots,\si_{i-1})
\int \pee(d\eta_{i+1}\cdots d\eta_{n} | \si_1,\ldots,\si_{i-1},\eta_i)\
g(\si_1,\ldots,\si_{i-1},\eta_i,\eta_{i+1},\ldots,\eta_{n})\leq
$$
$$
\max_{a\in A}
\int \pee(d\eta_{i+1}\cdots d\eta_{n} | \si_1,\ldots,\si_{i-1},\si_i=a)\
g(\si_1,\ldots,\si_{i-1},a,\eta_{i+1},\ldots,\eta_{n}) 
$$
$$
-\min_{b\in A}
\int \pee(d\eta_{i+1}\cdots d\eta_{n} | \si_1,\ldots,\si_{i-1},\si_i=b)\
g(\si_1,\ldots,\si_{i-1},b, \eta_{i+1},\ldots,\eta_{n})\,.
$$
\begin{equation}\label{torchon}
=: Y_i (\si) - X_i(\si)\,.
\end{equation}

Denote by $\pee^{\si}_{i,a,b}=
\pee^{\si_{{\scriptscriptstyle <i}}}_{i,a,b}$ the maximal coupling \cite{lindvall}
of the conditional distributions
$\pee(d\eta_{{\scriptscriptstyle \geq i+1}} | \si_1,\ldots,\si_{i-1},\si_i=a)$
and
$\pee(d\eta_{{\scriptscriptstyle \geq i+1}}| \si_1,\ldots,\si_{i-1},\si_i=b)$, that is,
the coupling for which 
Now we introduce the (infinite) upper-triangular matrix $D^\sigma$ defined for $i,j\in\N$ by
$$
D_{ii}^\si:=1
$$
\begin{eqnarray}\label{dracula}
D_{i,i+j}^\sigma
:=\max_{a,b\in A} \pee^{\si}_{i,a,b}\left\{\si^{(1)}_{i+j}\neq \si^{(2)}_{i+j}\right\}\,.
\end{eqnarray}
Notice that if the $\si_i$'s are mutually independent, then $D^\si$ is the identity matrix because
the conditional distributions 
$\pee(d\eta_{{\scriptscriptstyle \geq i+1}}| \si_1,\ldots,\si_{i-1},\si_i=a)$
and \\
$\pee(d\eta_{{\scriptscriptstyle \geq i+1}} | \si_1,\ldots,\si_{i-1},\si_i=b)$ are equal. Hence 
we have a perfect coupling in this case.

We proceed with the following simple telescoping identity:
$$
g(\si_1,\ldots,\si_{i-1},a,\si^{(1)}_{i+1},\ldots, \si^{(1)}_{n})
-g(\si_1,\ldots,\si_{i-1},b,\si^{(2)}_{i+1},\ldots,\si^{(2)}_{n})=
$$
$$
[g(\si_1,\ldots,\si_{i-1},a,\si^{(1)}_{i+1},\ldots,\si^{(1)}_{n})
-
g(\si_1,\ldots,\si_{i-1},b,\si^{(1)}_{i+1},\ldots,\si^{(1)}_{n})]+
$$
$$
[g(\si_1,\ldots,\si_{i-1},b,\si^{(1)}_{i+1},\ldots,\si^{(1)}_{n})-
g(\si_1,\ldots,\si_{i-1},b,\si^{(2)}_{i+1},\si^{(1)}_{i+2},\ldots,\si^{(1)}_{n})]+
$$
$$
[g(\si_1,\ldots,\si_{i-1},b,\si^{(2)}_{i+1},\si^{(1)}_{i+2},\ldots,\si^{(1)}_{n})-
g(\si_1,\ldots,\si_{i-1},b,\si^{(2)}_{i+1},\si^{(2)}_{i+2},\si^{(1)}_{i+3},\ldots,\si^{(1)}_{n})]
$$
$$
+\cdots +
$$
$$
[g(\si_1,\ldots,\si_{i-1},b,\si^{(2)}_{i+1},\si^{(2)}_{i+2},\ldots,\si^{(2)}_{n-1},\si^{(1)}_{n})-
g(\si_1,\ldots,\si_{i-1},b,\si^{(2)}_{i+1},\ldots,\si^{(2)}_{n})]
$$
$$
=:\sum_{j=0}^{n-i} \nabla^{12}_{i,i+j}g\,.
$$
We define the variation of $g$ at site $i$ by
$$
\delta_i g := \sup_{\substack{\si_j=\si'_j \\ \forall j\neq i}} |g(\si)-g(\si')|\,,
$$
and by  construction we have the inequality
$$
\nabla^{12}_{i,i+j}g \leq \delta_{i+j} g \ \1_{\si^{(1)}_{i+j}\neq \si^{(2)}_{i+j}}\,.
$$
It follows from \eqref{torchon} and \eqref{dracula} that
$$
Y_i(\si)-X_i(\si)=
$$
$$
\max_{a,b\in A}\Big\{
\int \pee(d\eta_{i+1}\cdots d\eta_{n} | \si_1,\ldots,\si_{i-1},\si_i=a)\
g(\si_1,\ldots,\si_{i-1},a,\eta_{i+1},\ldots,\eta_{n})
$$
$$
-\int \pee(d\eta_{i+1}\cdots d\eta_{n} | \si_1,\ldots,\si_{i-1},\si_i=b)\
g(\si_1,\ldots,\si_{i-1},b, \eta_{i+1},\ldots,\eta_{n})
\Big\}
$$
$$
=
\max_{a,b\in A}\Big\{
\int \pee^{\si}_{i,a,b}(d\si^{(1)}_{{\scriptscriptstyle \geq
i+1}},d\si^{(2)}_{{\scriptscriptstyle \geq i+1}})
$$
$$
\big[
g(\si_1,\ldots,\si_{i-1},a,\si^{(1)}_{i+1},\ldots,\si^{(1)}_{n})-
g(\si_1,\ldots,\si_{i-1},b, \si^{(2)}_{i+1},\ldots,\si^{(2)}_{n})
\big]\Big\}
$$
$$
\leq \max_{a,b\in A} \sum_{j=0}^{n-i} \delta_{i+j}g\ 
\pee^{\si}_{i,a,b}\left\{\si^{(1)}_{i+j}\neq \si^{(2)}_{i+j}\right\}
$$
$$
\leq \sum_{j=0}^{n-i} D_{i,i+j}^\si\ \delta_{i,i+j}g=(D^{\si}\delta g)_i
$$
where $\delta g$ denotes the column vector with coordinates $\delta_j g$, for
$j=1,\ldots,n$, and $0$ for $j>n$.
Therefore, we get the inequality
\begin{equation}\label{charlo}
V_i(\si)=Y_i(\si) - X_i(\si) \leq (D^\si \delta g)_i\,.
\end{equation}
Applying the above reasoning to $-g$ shows that the previous
inequality also applies to $-V_i$.

\begin{remark}
The advantage of the previous bound is that it only involves $\delta g$.
One could imagine to consider, for instance, the second moment of
$\nabla^{12}_{i,i+j} g$ instead. This could lead to better results
but it has the drawback that we need to know much more about the coupling
than we usually do.
\end{remark}

\subsection{Uniform decay of $D^\si$: exponential concentration inequality}

Let $\overline{D}_{i,j}:= \sup_{\si\in A^{\N}} D_{i,j}^\sigma$. We assume that
the following operator $\ell^2(\N)$-norm is finite:
\begin{equation}\label{dobru}
\|\overline{D} \|^2_{\ldeux} := \sup_{u\in \ell_2(\N),
\|u\|_{\ldeux}=1} \| \overline{D} u\|^2_{\ldeux} <\infty\,.
\end{equation}

We have the following exponential concentration inequality.

\begin{theorem}\label{expconthm}
Let $n\in\N$ be arbitrary.
Assume that \eqref{dobru} holds.
Then, for all functions $g:A^n\to\R$,
we have the inequality
\begin{equation}\label{expcon}
\pee\left\{|g - \E g| \geq t \right\} \leq
2 \exp\left(-\frac{2 t^2}{\| \overline{D} \|^2_{\ldeux} \ \|\delta g \|_{\ldeux}^2}\right)\cdot
\end{equation}
for all $t>0$.
\end{theorem}

\begin{proof}
We recall the following lemma which is proved  in \cite{dl}.
\begin{lemma}
Suppose $\mathcal{F}$ is a sigma-field and $Z_1,Z_2,V$ are random
variables such that
\begin{enumerate}
\item $Z_1\leq V\leq Z_2$
\item $\E (V|\mathcal{F}) =0$
\item $Z_1$ and $Z_2$ are $\mathcal{F}$-measurable.
\end{enumerate}
Then, for all $\lambda\in\R$,  we have
\begin{equation}\label{kwats}
\E( e^{\lambda V}|\mathcal{F}) \leq e^{\lambda^2(Z_2-Z_1)^2/8}\,.
\end{equation}
\end{lemma}
We apply this lemma with $V=V_i$, $\mathcal{F}=\mathcal{F}_{i-1}$,

$Z_1=X_i- \E[g|\mathcal{F}_{i-1}]$, $Z_2=Y_i-
\E[g|\mathcal{F}_{i-1}]$. Using inequality \eqref{charlo}
$$
V_i(\si)=Y_i(\si) -X_i(\si) \leq (D^\si \delta g)_i
$$
we obtain
\begin{equation}
\E (e^{\lambda V_i}|\mathcal{F}_{i-1})(\si)\leq e^{\lambda^2 (D^\si \delta g)_i^2/8}\,.
\end{equation}
Therefore, by successive conditioning, and the exponential Chebychev's inequality,
\begin{eqnarray}
\pee\left\{ g -\E g \geq t\right\}
& \leq &
e^{-\lambda t} \E\left(e^{\lambda \sum_{i=1}^n V_i}\right)
\nonumber
\\
&\leq &
e^{-\lambda t} \E\left(\E(e^{\lambda V_n}|\mathcal{F}_{n-1}) e^{\lambda\sum_{i=1}^{n-1}V_i}\right)
\nonumber\\
&\leq &\!\!\!\!\cdots \leq
e^{-\lambda t} \exp\left( \frac{\lambda^2}{8} \| \overline{D} \delta
g\|^2_{\ldeux} \right)\\
\nonumber
& \leq &
e^{-\lambda t} \exp\left( \frac{\lambda^2}{8} \| \overline{D}\|^2_{\ldeux}\
\|\delta g\|^2_{\ldeux}\right). 
\end{eqnarray}
Now choose the optimal $\lambda=4t/(\| \overline{D}\|^2_{\ldeux} \ \|\delta
g\|^2_2)$ to obtain
$$
\pee\left\{g - \E g \geq t \right\} \leq
\exp\left(-\frac{2 t^2}{\| \overline{D} \|_{\ldeux}^2 \ \|\delta g \|_{\ldeux}^2}\right)\cdot
$$
Combining the inequality for $g$ and the one for $-g$ yields \eqref{expcon}. The theorem is proved.
\end{proof}

\subsection{Non-uniform decay of $D^\si$: polynomial and
stretched-exponential concentration inequalities }

If the dependence on $\si$ of the elements of the coupling matrix 
cannot be controlled uniformly, then in many cases we can still control
the moments of the coupling matrix.
To this aim, we introduce the (non-random, i.e., not depending on $\si$) matrices
\begin{equation}\label{drond}
\D^{(p)}_{i,j}:= \E[(D_{i,j}^\sigma)^p]^{1/p}
\end{equation}
for all $p\in\N$.

A typical example of non-uniformity which we will encounter, for instance in the
low-temperature Ising model, is an estimate of the following form:
\begin{equation}\label{li}
D_{i,i+j}^\si \leq \1 \{\ell_i(\si)\geq j\} + \psi_{j}
\end{equation}
where $\psi_{j}\ge 0$ does not depend on $\si$, and where $\ell_i$ are
unbounded functions of $\si$ with a distribution independent of $i$.
The idea is that the matrix elements $D_{i,i+j}^\si$ ``start to decay'' when $j\geq \ell_i(\si)$. The ``good''
configurations $\si$ are those for which $\ell_i(\si)$ is ``small''.

In the particular case when \eqref{li} holds, in principle one still {\em can} have
an exponential concentration inequality provided one is able to bound
$$
\E\left(
e^{\lambda \sum_{i=1}^{n} \ell_i^2}\right)\,.
$$
However, in the example given below, the tail of the $\ell_i$ will be
stretched exponential. Henceforth, we cannot deduce an
exponential concentration inequality from these estimates.

We now prove an inequality for the variance of $g$
which is a generalization of an inequality derived in \cite{devroye}
in the i.i.d. case.

\begin{theorem}\label{devroyeineq}
Let $n\in\N$ be arbitrary. Then for all functions
$g:A^n\to\R$ we have the inequality
\begin{equation}\label{DEV}
\E\left[ 
(g -\E g )^2 \right]
\leq
\| \D^{(2)}\|^2_{\ldeux}\ \| \delta g\|_{\ldeux}^2\,.
\end{equation}
\end{theorem}

\begin{proof}
We start again from the decomposition \eqref{decomp}. 
Recall the fact that $\E[V_i |\mathcal{F}_j]=0$ for all $i> j$, from which it follows
that $\E[V_i V_j]=0$ for $i\neq j$.
Using \eqref{charlo} and Cauchy-Schwarz's inequality we obtain
\begin{eqnarray}
\nonumber
&\E& \!\!\!\!\left[ (g -\E g)^2\right]
=  \E \sum_{i=1}^n V_i^2\\
\nonumber
& \leq & \E\left(\sum_{i=1}^n (D \delta g)_i^2\right)\\
\nonumber
&= & \sum_{i=1}^n \sum_{k=1}^n \sum_{l=1}^n \E\left(D_{i,k} D_{i,l}\right) \delta_k g
\delta_l g \\
\nonumber
&\leq & \sum_{i=1}^n \sum_{k=1}^n \sum_{l=1}^n \E\left( D_{i,k}^2\right)^{\frac{1}{2}}
\E\left( D_{i,l}^2\right)^{\frac{1}{2}}\delta_k g \delta_l g \\
\nonumber
& = & \| \D^{(2)}\delta g\|_{\ldeux}^2\\
\nonumber
&\leq & \| \D^{(2)}\|^2_{\ldeux}\ \| \delta g\|_{\ldeux}^2\,.
\end{eqnarray}
\end{proof}

\begin{remark}
In the i.i.d. case, the coupling matrix $D$ is the identity
matrix. Hence inequality \eqref{DEV} reduces to
$$
\E\left[
(g -\E g)^2 \right]
\leq
\|\delta g\|_{\ldeux}^2
$$
which is the analogue of Theorem 4 in \cite{devroye}.
\end{remark}

We now turn to higher moment estimates. We have the following theorem
from which we recover Theorem \ref{devroyeineq} but with a bigger
constant.

\begin{theorem}\label{moments}
Let $n\in\N$ be arbitrary. For all functions $g:A^n\to\R$ and
for any $p\in\N$, we have
$$
\E\left[ (g-\E g)^{2p}\right]\leq
(20p)^{2p} \ \|\D^{(2p)}\|^{2p}_{\ldeux} \|\delta g\|^{2p}_{\ldeux}\,.
$$
\end{theorem}

\begin{proof}
We start from \eqref{decomp} and get
$$
\E\left[ (g-\E g)^{2p}\right]=
\E\left[\left(\sum_{i=1}^n V_i\right)^{2p}\right]\,.
$$
Now, by \eqref{decomp} and since $\E(V_{i}|\mathcal{F}_{j})=0$ for $i>j$, 
$g-\E g$ is a martingale, to which we apply Burkholder-Gundy's
inequality \cite[formula II.2.8, p. 41]{garsia}: for any $q\ge 2$, we
have
$$
\E\left[|g-\E g|^q\right]^{\frac{1}{q}}\le 10 q \
\E\left[\left(\sum_{i=1}^{n}V_{i}^{2}\right)^{\frac{q}{2}}\right]^{\frac{1}{q}}\;.
$$
Therefore, for $q=2p$, $p\in\N$, this gives at once
$$
\E\left[ (g-\E g)^{2p}\right]\leq (20 p)^{2p} \ \E\left[\left(\sum_{i=1}^n
V_i^2\right)^{p}\right]\, .
$$
We now estimate the rhs by using \eqref{charlo}:
\begin{eqnarray}
& \E & \!\!\!\! \big[(\sum_{i}V_i^2)^p\big]\\
& = & \sum_{i_1}\cdots \sum_{i_{p}} \E\left( V_{i_1}^2\cdots
  V_{i_{p}}^2\right)\\
\nonumber 
& \leq & 
\sum_{i_1}\cdots \sum_{i_{p}} \E\left[ (D \delta g)_{i_1}^2\cdots
  (D \delta g)_{i_{p}}^2\right]\\
\nonumber 
& =& 
\sum_{i_1\cdots i_p} \sum_{j_1\cdots j_p} \sum_{k_1\cdots k_p}
\E\left(
\prod_{r=1}^{p} D_{i_r,j_r} D_{i_r,k_r}\right)\
\left(\prod_{r=1}^{p} \delta_{j_r}g \ \delta_{k_r}g \right)\\
\nonumber 
& \leq & 
\sum_{i_1\cdots i_p} \sum_{j_1\cdots j_p} \sum_{k_1\cdots k_p}
\prod_{r=1}^{p}\left(
\D^{(2p)}_{i_r,j_r} \D^{(2p)}_{i_r,k_r} \delta_{j_r} g\ \delta_{k_r} g
\right)\\
\nonumber 
& =& 
\| \D^{(2p)}\delta g\|^{2p}_{\ldeux} \leq \| \D^{(2p)}\|^{2p}_{\ldeux} \ \| \delta g\|^{2p}_{\ldeux}
\end{eqnarray}
where in the fourth step we used the inequality
$$
\E(f_1 \cdots f_{2p}) \leq \prod_{i=1}^{2p} (\E(f_i^{2p}))^{\frac{1}{2p}}
$$
which follows from H\"older's inequality.
\end{proof}

In order to be able to apply Theorems \ref{devroyeineq} and \ref{moments}, one needs to
estimate $\|\mathcal{D}^{(2p)}\|_{\ldeux}$. 

\begin{proposition}\label{prop2} Assume inequality (\ref{li}) holds, and
let $p\in\N$. We have the bound
$$
\|\mathcal{D}^{(2p)}\|_{\ldeux}\le \sum_{j=1}^\infty\pee\big(
\ell_{0}(\sigma)\ge   j\big)^{1/2p}+\|\psi\|_{\lun} \;. 
$$
\end{proposition}

\begin{proof}
We start by an upper estimate of $\mathcal{D}^{(2p)}$. From the
definition \eqref{drond} and the bound (\ref{li}) we have using Minkowski's  inequality (for $j\ge i$)
$$
\mathcal{D}^{(2p)}_{i,j}=\E\big[(D_{i,j}^\sigma)^{2p}\big]^{1/2p}\le 
\E\Big(\big(\1\big\{ \ell_{i}(\sigma)\ge
  j-i\big\}+\psi_{j-i}\big)^{2p}\Big)^{\frac{1}{2p}}\le
$$
\begin{equation}\label{merlot}
\E\left[\big(\1\big\{ \ell_{i}(\sigma)\ge
  j-i\big\}\big)^{2p}\right]^{\frac{1}{2p}}+\psi_{j-i}\le
\pee\big( \ell_{0}(\sigma)\ge
j-i\big)^{\frac{1}{2p}}+\psi_{j-i}=:u_{i-j}\;,
\end{equation}
since the law of $ \ell_{i}$ is independent of $i$.

Now take $v\in\ell^2(\N)$ with $\|v\|_{\ldeux}=1$.
We have
$$
\| \mathcal{D}^{(2p)}v\|_{\ldeux} \leq 
\left\| \sum_{k=1}^\infty \mathcal{D}^{(2p)}_{i,k} |v_k|\right\|_{\ldeux}
\leq \left\| \sum_{k=1}^\infty u(i-k) |v_k|\right\|_{\ldeux}
$$
where the second inequality comes from \eqref{merlot}. Since
we have the $\ell^2(\N)$-norm of a convolution, we can apply
Young's inequality (see, e.g., \cite{zyg}) to get
$$
\|\mathcal{D}^{(2p)}\|_{\ldeux} \leq \|u\|_{\lun}\,.
$$ 
The result immediately follows.
\end{proof}

Before we state the next theorem, which is a corollary of 
Proposition \ref{prop2} and Theorem
\ref{moments}, we need the definition of some Orlicz spaces.
We only deal here with a restricted class useful in our applications, 
we refer to  \cite{raoren,zyg} for the
general definition. For $\varrho>0$, let $\Phi_{\varrho}:\R\to\R^+$ be the Young
function defined by
$$
\Phi_{\varrho}(x)=e^{(|x|+h_{\varrho})^{\varrho}}-
e^{h_{\varrho}^{\varrho}}
$$
where $h_{\varrho}=((1-\varrho)/\varrho)^{1/\varrho}\1\{0<\varrho<1\}$. These are the Young
functions used in particular in \cite{su}.
We recall that (see \cite{zyg})
the Luxembourg norm with respect to $\Phi_{\varrho}$ of a
random variable $Z$ is defined by
$$
\big\| Z\big\|_{\Phi_{\varrho}}=\inf\left\{\lambda>0\;\bigg|\;
\E\left(\Phi_{\varrho}\left(\frac{Z}{\lambda}\right)\right)\le
1\right\}\;.
$$

\begin{remark}
Note that for $\Phi_p(x)=|x|^p$, the Luxembourg norm is nothing
but the $L^p$-norm.
\end{remark}

\begin{theorem}\label{curanto}
Let $n\in\N$ be arbitrary. Then, for all functions $g:A^n\to\R$,
for any $p\in\N$ and any $\epsilon>0$, we have
$$
\E\big[(g-\E g)^{2p}\big]\le
$$
\begin{equation}\label{harvest}
(20 p)^{2p}\left(\zeta(1+\epsilon/(2p-1))^{(2p-1)/2p}
\E\big(\ell_{0}^{2p+\epsilon}\big)^{1/2p}
+\|\psi\|_{\lun}\right)^{2p}\|\delta g\|_{\ldeux}^{2p}
\end{equation}
where, as usual, $\zeta$ denotes Riemann's zeta function.
For any $\vartheta>0$, there is a constant $C_{\vartheta}>0$,  such that 
for any $\varrho<\vartheta/(1+\vartheta)$ satisfying
$\zeta(\vartheta(1-\varrho)/\varrho)\ge 1$, we have
\begin{equation}\label{cabernet}
\left\|\frac{g-\E g}{\;\|\delta g\|_{\ldeux}}
\right\|_{\Phi_{\varrho}}
\le C_{\vartheta}\left(
\zeta\left(\frac{\vartheta(1-\varrho)}{\varrho}\right)
\|\ell_{0}\|_{\Phi_{\vartheta}}^{\vartheta(1-\varrho)/\varrho}
+\|\psi\|_{\lun}\right)\;.
\end{equation}
\end{theorem}

\begin{remark}
A similar result holds when $\zeta(\vartheta(1-\varrho)/\varrho)< 1$
with the square root of the zeta function. 
Note also that when $\varrho$ increases to
$\vartheta/(1+\vartheta)$, the number $\vartheta(1-\varrho)/\varrho$
decreases to one.  
\end{remark}

\begin{proof}
We first estimate $\|\mathcal{D}^{(2p)}\|_{\ldeux}$ in terms of some
moment of $\ell_{0}$. Let $\epsilon'=\epsilon/(2p-1)$. We have using
H\"older inequality
\begin{eqnarray*}
\sum_{j=1}^\infty \pee\big(\ell_{0}(\sigma) \ge j\big)^{1/2p} & = &
\sum_{j=1}^\infty j^{(2p-1)(1+\epsilon')/2p}
\pee\big(\ell_{0}(\sigma)\ge   j\big)^{1/2p}\;
j^{-(2p-1)(1+\epsilon')/2p}\\
& \le &\zeta(1+\epsilon')^{(2p-1)/2p}\left(\sum_{j=1}^\infty j^{2p-1+\epsilon}
\pee\big(\ell_{0}(\sigma)\ge   j\big)\right)^{1/2p}\\
& \le &\zeta(1+\epsilon')^{(2p-1)/2p}\
\E\big(\ell_{0}^{2p+\epsilon}\big)^{1/2p} \;.
\end{eqnarray*}
Using Proposition \ref{prop2} we get 
$$
\|\mathcal{D}^{(2p)}\|_{\ldeux}
\le  \zeta(1+\epsilon')^{(2p-1)/2p}\
\E\big(\ell_{0}^{2p+\epsilon}\big)^{1/2p}
+\|\psi\|_{\lun}\;,
$$
and \eqref{harvest} follows using Theorem
\ref{moments}. 

To prove \eqref{cabernet}, we first observe that from \eqref{harvest}
we have, for  $q$ even 
$$
\left\|\frac{g-\E g}{\;\|\delta g\|_{\ldeux}}\right\|_{\Lq}
\le 10 q  \Big(\zeta(1+\epsilon/(q-1))^{\frac{q-1}{q}}
\E(\ell_{0}^{q+\epsilon})^{\frac{1}{q}}+
\|\psi\|_{\lun}\Big)\;.
$$
We now  recall that for any $1>\varrho>0$, there
is a constant $\tilde C_{\varrho}>1$ such that 
$$
\tilde C_{\varrho}^{-1}\sup_{q>2}\frac{\|Z\|_{\Lq}}{q^{1/\varrho}}
\le \|Z\|_{\Phi_{\varrho}}\le \tilde 
C_{\varrho}\sup_{q>2}\frac{\|Z\|_{\Lq}}{q^{1/\varrho}}\cdot
$$
(See, e.g.,  \cite{su} for a proof.)  
It is easy to verify  using Young's inequality that the same
inequality holds (with slightly different constants)
when the supremum is taken over the $q$ even
integers, and we will only consider such $q$ below. 

Therefore if $0<\varrho< \vartheta/(1+\vartheta)$,  taking
$$
\epsilon=\vartheta q \left(\frac{1}{\varrho}-\frac{1}{\vartheta}-1\right)
$$
we get 
\begin{eqnarray*}
\left\|\frac{g-\E g}{\;\|\delta g\|_{\ldeux}}
\right\|_{\Phi_{\varrho}} & \le & \Oun 
\sup_{q>2} q^{1-1/\varrho}
\zeta(1+\epsilon/(q-1))^{(q-1)/q}\
\E(\ell_{0}^{q+\epsilon})^{1/q}
+\Oun \|\psi\|_{\lun}\\
& \le & \Oun 
\sup_{q>2} 
\zeta\left(1+\frac{q(\vartheta-\varrho-\vartheta\varrho)}{\varrho(q-1)}\right)^{(q-1)/q}
\|\ell_{0}\|_{\Phi_{\vartheta}}^{\vartheta(1-\varrho)/\varrho}
+\Oun \|\psi\|_{\lun}\\
& \le & \Oun 
\sup_{q>2} 
\zeta\left(\frac{\vartheta(1-\varrho)}{\varrho}\right)^{(q-1)/q}
\|\ell_{0}\|_{\Phi_{\vartheta}}^{\vartheta(1-\varrho)/\varrho}
+\Oun \|\psi\|_{\lun}
\end{eqnarray*}
since the function $\zeta$ is decreasing. Thus \eqref{cabernet} is proved.
The proof of the theorem is now complete.
\end{proof}

It is easy to obtain from Theorem \ref{curanto} the following
concentration inequalities.

\begin{proposition}\label{pinard}
Let $n$ be an arbitrary positive integer. 

\begin{itemize}
\item If $\E(\ell_0^{2p+\epsilon})<\infty$ (for some $\epsilon>0$), and
$\|\psi\|_{\lun}<\infty$, we have
\begin{equation}\label{moule}
\pee\left\{ |g -\E g| >t\right\} \leq C_p \frac{\| \delta
g\|_{\ldeux}^{2p}}{t^{2p}}
\end{equation}
where $C_p\in ]0,\infty[$, $p\in\N$, for any $g:A^n\to\R$.

\item
Let $0<\varrho<1$. If $\|\psi\|_{\lun}<\infty$, and 
$\|\ell_{0}\|_{\Phi_{\vartheta}}<\infty$
for some $\vartheta>\varrho/(1-\varrho)$,
there exists a constant $c_{\varrho,\vartheta} \in]0,\infty[$
such that 
\begin{equation}\label{huitre}
\pee\left\{ |g -\E g| >t\right\} \leq 4
\exp\left(-c_{\varrho,\vartheta}\ \frac{t^\varrho}{\|\delta
  g\|_{\ldeux}^\varrho}\right)\;,
\end{equation}
for any $g:A^n\to\R$.

\end{itemize}
\end{proposition}

\begin{proof}
The proof of \eqref{moule} is an immediate consequence of \eqref{harvest} applied
to $g$ and $-g$ and Chebychev's inequality. 

For the proof of \eqref{huitre}, we have for any $\lambda>0$
 using Chebychev's inequality
\begin{eqnarray*}
\pee( g -\E g > t) & = &
\pee\left( \frac{g -\E g}{\lambda \|\delta g\|_{\ldeux}} 
> \frac{t}{\lambda \|\delta g\|_{\ldeux}}\right)\\ 
& \leq & \pee
\left( \Phi_\varrho\left( 
\frac{g- \E g}{\lambda\|\delta g\|_{\ldeux}  }\right) >
\Phi_\varrho\left(\frac{t}{\lambda \|\delta g\|_{\ldeux}}\right)\right)\\ 
& \leq  &\frac{1}{\Phi_\varrho \big(t/(\lambda \|\delta g\|_{\ldeux})\big)}\
\E\left[ \Phi_\varrho\left( \frac{g- \E g}{\lambda 
\|\delta g\|_{\ldeux}}\right)\right]\;.
\end{eqnarray*}
We now take $\lambda=\|(g -\E
g)/\|\delta g\|_{\ldeux} \|_{\Phi_\varrho}$. 
By definition, $\E\left[\Phi_\varrho
\left(\frac{g-\E g}{\lambda\|\delta g\|_{\ldeux} }\right)\right]=1$. 
Thus we have
$$
\pee( g -\E g > t)\leq \frac{1}{\Phi_\varrho (t/\|g -\E g\|_{\Phi_\varrho})}\cdot
$$
Of course, the same inequality holds with $-g$. Applying
\eqref{cabernet} yields \eqref{huitre}. 
The proposition is proved. 
\end{proof}

In concrete applications of inequality \eqref{moule} we have to check
that $C_p<\infty$, otherwise the inequality is useless. To apply
\eqref{huitre}, we have to check that $c_\varrho>0$. We will give
an example of application below. 

Inequality \eqref{moule} is a ``polynomial'' concentration inequality
whereas inequality \eqref{huitre} is a ``stretched-exponential'' concentration inequality.

\begin{remark}
The 4 in the r.h.s. of \eqref{huitre} is not optimal. It can be
replaced by $2/(1-\epsilon)$ for any $\epsilon\in]0,1[$.
\end{remark}

\subsection{Random fields}\label{RF}

We now present the extension of our previous results to random fields.
This requires mainly notational changes.
We work with lattice spin systems. The configuration space is
$\Omega = \{ -, + \}^{\Zd}$, endowed with the product topology. 
We could of course take any finite set $A$ instead of $\{-,+\}$.
For $\la\subset\Zd$ and $\si,\eta\in\Omega$ we denote
$\si_\la\eta_{\la^c}$ the configuration coinciding with
$\si$ (resp.\ $\eta$) on $\la$ (resp.\ $\la^c$). 
A local function $g:\Omega\to\R$ is such that there
exists a finite subset $\Lambda\subset\Zd$ such that
for all $\si,\eta,\omega$, $g(\si_\Lambda\omega_{\Lambda^c}) = g(\si_\Lambda\eta_{\Lambda^c})$.

For $\si\in\Omega$ and $x\in\Zd$, $\si^x$ denotes the configuration
obtained from $\si$ by ``flipping'' the spin at $x$.
We denote $\delta_x g = \sup_\si |g(\si^x) - g(\si)|$ the variation of
$g$ at $x$. $\delta g$ denotes the map $\Zd\to\R: x\mapsto \delta_x g$.

We introduce the spiraling enumeration $\Gamma:\Zd\to\N$ illustrated in the
figure for the case $d=2$.

\begin{figure}[!h]
\begin{center}
\epsfig{width=7cm,clip=,figure=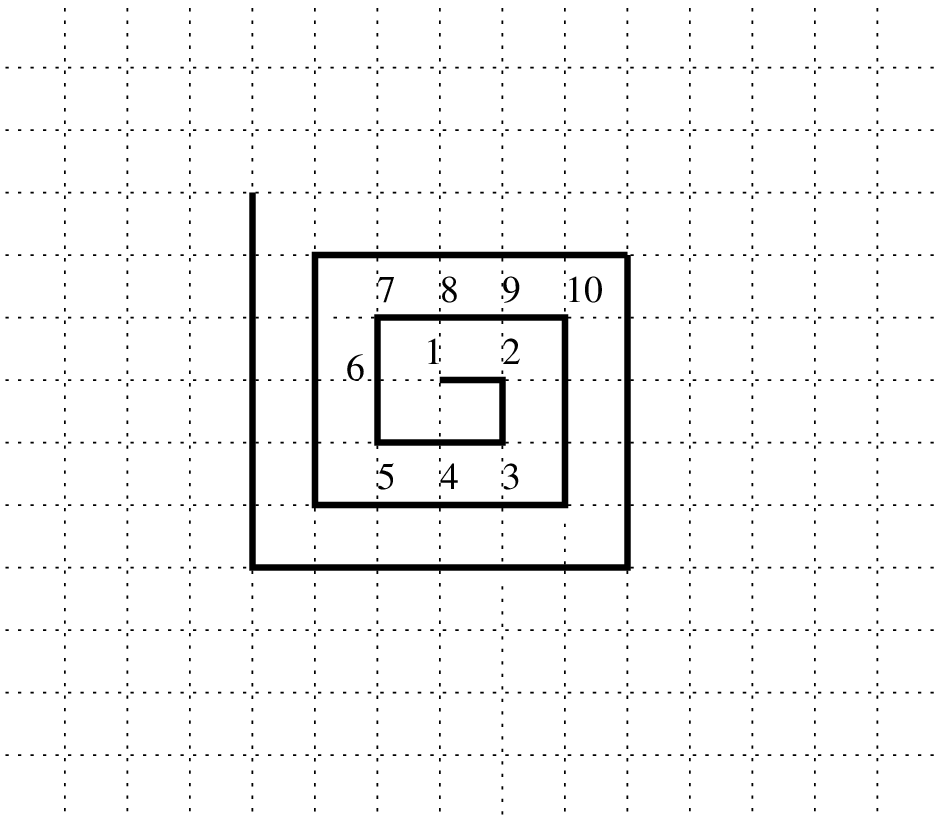}
\end{center}
\end{figure}

We will use the abbreviation $(\leq x)=\{y\in\Zd : \Gamma(y)\leq \Gamma(x)\}$ and similarly we introduce
the abbreviations $(<x)$.
By definition $\mathcal{F}_{\leq x}$ denotes the sigma-field generated by $\si(y)$, $y\leq x$ and
$\mathcal{F}_{<0}$ denotes the trivial sigma-field. 

For any local function $g:\Omega\to\R$, we have the analog decomposition as in \eqref{decomp}:
\begin{equation}\label{osition}
g - \E g= \sum_{x\in\Zd} V_x
\end{equation}
where 
$$
V_x:= \E\left[g|\mathcal{F}_{\leq x}]- \E[g|\mathcal{F}_{<x}\right]\,.
$$
The analog of the coupling matrix is the following matrix indexed by lattice sites $x,y\in\Zd$
\begin{equation}\label{cigare}
D_{x,y}^\si:= \hat{\pee}_{x,+,-}^{\si}\left\{X_{1}(y)\neq X_{2}(y)\right\}
\end{equation}
where
$
\hat{\pee}_{x,+,-}^{\si}
$
denotes the maximal coupling between the conditional measures
$\pee(\cdot | \si_{<x,+_x})$ and $\pee(\cdot | \si_{<x,-_x})$. The
notation ``$+_x$'' (resp. ``$-_x$'') means that at coordinate $x$ 
in the configuration we put a ``$+$'' (resp. a ``$-$'').

We first consider the case of uniform decay of $D$.
In that case, the exponential concentration inequality of Theorem \ref{expconthm}
holds with the norm of $\ell_2(\Zd)$, i.e.,
$\|\delta g\|_2^2= \sum_{x\in\Zd} (\delta_x g)^2$ (which is trivially
finite since $g$ is a local function).

\begin{theorem}\label{chardon}
Assume that 
\begin{equation}\label{the}
\overline{D}_{x,y}:= \sup_{\si} D_{x,y}^\si
\end{equation}
is a bounded operator in $\ell_2(\Zd)$.
Then for all local functions $g$ we have the following inequality
\begin{equation}\label{expconzd}
\pee\left\{|g - \E g| \geq t \right\} \leq 2
\exp\left(-\frac{2 t^2}{\|\overline{D} \|_{\ldeuxzd}^2 \ \|\delta g \|_{\ldeuxzd}^2}\right)
\end{equation}
for all $t>0$.
\end{theorem}

In the non-uniform case, Theorems \ref{moments},
\ref{curanto} and Proposition \ref{pinard}
extend immediately as follows. The analog of \eqref{li} is 
\begin{equation}\label{lizd}
D_{x,y}^\si \leq \1\{\ \ell_x(\si) \geq |y-x| \} + \psi(|y-x|)\;.
\end{equation}
From now on, we  assume
that the distribution of $\ell_x$ 
is independent of $x$. We extend the matrix $\mathcal{D}$ defined in \eqref{drond} by putting
$$
\D^{(p)}_{x,y}:= \E[(D_{x,y}^\sigma)^p]^{1/p}
$$
for $x,y\in\Zd$.

\begin{theorem}\label{moments_bis}
For any local function $g$ and
for any $p\in\N$, we have
$$
\E\left[ (g-\E g)^{2p}\right]\leq
(20p)^p \ \|\D^{(2p)}\|^{2p}_{\ldeuxzd} \|\delta g\|^{2p}_{\ldeuxzd}\,.
$$
\end{theorem}

\begin{theorem}\label{curantobis}
For any local function $g$,
for any $p\in\N$ and any $\epsilon>0$, we have
$$
\E\big[(g-\E g)^{2p}\big]\le
$$
\begin{equation}\label{harvestbis}
(20 p)^{2p}\left(\zeta(1+\epsilon/(2p-1))^{(2p-1)/2p}
\E\big(\ell_{0}^{2pd+\epsilon}\big)^{1/2p}
+\|\psi\|_{\lun}\right)^{2p}\|\delta g\|_{\ldeuxzd}^{2p}
\end{equation}
where, as usual, $\zeta$ denotes Riemann's zeta function.
For any $\vartheta>0$, there is a constant $C_{\vartheta}>0$,  such that 
for any $\varrho<\vartheta/(1+\vartheta)$ satisfying
$\zeta(\vartheta(1-\varrho)/\varrho)\ge 1$, we have
\begin{equation}\label{cabernetbis}
\left\|\frac{g-\E g}{\;\|\delta g\|_{\ldeuxzd}}
\right\|_{\Phi_{\varrho}}
\le C_{\vartheta}\left(
\zeta\left(\frac{\vartheta(1-\varrho)}{\varrho}\right)
\|\ell_{0}^d\|_{\Phi_{\vartheta}}^{\vartheta(1-\varrho)/\varrho}
+\|\psi\|_{\lun}\right)\;.
\end{equation}
\end{theorem}

\begin{proposition}\label{pinardbis}
For any local function $g$ we have the inequalities:

\begin{itemize}
\item If $\E(\ell_0^{2pd+\epsilon})<\infty$ (for some
$\epsilon>0$), and $\|\psi\|_{\lun}<\infty$, we have
$$
\pee\left\{ |g -\E g| >t\right\} \leq C_p \frac{\| \delta
g\|_{\ldeuxzd}^{2p}}{t^{2p}}
$$
where $C_p\in ]0,\infty[$, $p\in\N$. 

\item
Let $0<\varrho<1$. If $\|\psi\|_{\lun}<\infty$, and 
$\|\ell_{0}^d\|_{\Phi_{\vartheta}}<\infty$
for some $\vartheta>\varrho/(1-\varrho)$,
there exists a constant $c_{\varrho,\vartheta} \in]0,\infty[$
such that 
$$
\pee\left\{ |g -\E g| >t\right\} \leq 4
\exp\left(-c_{\varrho,\vartheta}\ \frac{t^\varrho}{\|\delta g\|_{\ldeuxzd}^\varrho}\right)\cdot
$$
\end{itemize}
\end{proposition}

\begin{remark}
It is immediate to extend the previous inequalities to integrable functions $g$ belonging to the
closure of the set of local functions with the norm $|\!|\!| g |\!|\!|:=\|\delta g\|_{\ldeuxzd}$.
\end{remark}

\subsection{Existence of the coupling by bounding the variation}

We continue with random fields and state a proposition which says that if we have an
estimate of the form
$$
V_x \leq (D\delta g)_x
$$
for some matrix $D$, then there exists a coupling with coupling matrix $\hat{D}$ such that
its matrix elements decay at least as fast as the matrix elements of $D$. 
We formulate the proposition more abstractly:

\begin{proposition}\label{sardine}
Suppose that $\pee$ and $\qee$ are probability measures on $\Omega$
and $g:\Omega\to\mathbb{R}$ such that we have the
estimate
\begin{equation}\label{extasy}
\left|
\E_\pee[g]-\E_\qee[g]
\right|
\leq 
\sum_{x\in\Zd} \rho(x) \delta_x g
\end{equation}
for some ``weights'' $\rho:\Zd\to \R^+$.
Suppose $\varphi:\Zd\to \R^+$ is such that 
$$
\sum_{x\in\Zd} \rho(x) \varphi(x) <\infty \,.
$$
Then there exists a coupling $\hat{\mu}$ of $\pee$ and $\qee$ such that 
$$
\sum_{x\in\Zd} 
\hat{\mu} \left\{ X_1(x)\neq X_2(x)\right\} \varphi(x) \leq
\sum_{x\in\Zd} \varphi(x)\rho(x)<\infty\,.
$$
\end{proposition}

\begin{proof}
Let $B_n:=[-n,n]^d \cap \Zd$. Define the ``cost'' function
$$
C_n^\varphi(\si,\si'):= \sum_{x\in B_n} |\si_x -\si'_x|\ \varphi(x)\,.
$$
Denote by $\pee_n$, resp. $\qee_n$, the joint distribution of $\{\si_x, x\in B_n\}$
under $\pee$, resp. $\qee$.
Consider the class of functions
$$
\mathcal{G}_{C_n^\varphi} :=\{
g| \ g\in \mathcal{F}_{B_n},\, |g(\si)-g(\si')| \leq \sum_{x\in\Zd} \varphi(x) \1\{\si_x\neq \si'_x\},\,
\forall \si,\si'\in \Omega\}\,.
$$
It is obvious from the definition that $g\in\mathcal{G}_{C_n^\varphi}$, if, and only if,
$g$ is $\mathcal{F}_{B_n}$-measurable and
$$
(\delta_x g)(\si)\leq \varphi(x)\quad \forall x\in B_n,\, \forall \si\in \Omega\,.
$$
Therefore, if \eqref{extasy} holds, then for all $g\in\mathcal{G}_{C_n^\varphi}$, 
$$
\left|
\E_\pee[g]-\E_\qee[g]
\right|
\leq 
\sum_{x\in\Zd} \rho(x) \delta_x g \leq \sum_{x\in\Zd} \rho(x) \varphi(x) \,. 
$$
Hence, by the Kantorovich-Rubinstein duality theorem \cite{rachev},
there exists a coupling $\hat{\mu}_n$ of $\pee_n$ and $\qee_n$ such that
$$
\E_{\hat{\mu}_n}\left(C_n^\varphi(\si,\si')\right)=
\E_{\hat{\mu}_n}\left(\sum_{x\in B_n} \varphi(x)\1\{X_1(x)\neq
X_2(x)\}\right)\leq \sum_{x\in\Zd} \varphi(x) \rho(x)\,. 
$$
By compactness (in the weak topology), there exists a subsequence along which $\hat{\mu}_n$
converges weakly to some probability measure $\hat{\mu}$. For any $k\leq n$, we have
$$
\E_{\hat{\mu}_n}\left(\sum_{x\in B_k} \varphi(x)\1\{X_1(x)\neq X_2(x)\}\right)\leq
$$
$$
\E_{\hat{\mu}_n}\left(\sum_{x\in B_n} \varphi(x)\1\{X_1(x)\neq X_2(x)\}\right)\leq 
\sum_{x\in\Zd} \varphi(x) \rho(x)\,.
$$
Therefore, taking the limit $n\to\infty$ along the above subsequence yields
$$
\E_{\hat{\mu}}\left(\sum_{x\in B_k} \varphi(x)\1\{X_1(x)\neq X_2(x)\}\right)\leq 
\sum_{x\in\Zd} \varphi(x) \rho(x)\,.
$$
We now take the limit $k\to\infty$ and use monotonicity to conclude that
$$
\E_{\hat{\mu}}\left(\sum_{x\in \Zd} \varphi(x)\1\{X_1(x)\neq X_2(x)\}\right)\leq 
\sum_{x\in\Zd} \varphi(x) \rho(x)\,.
$$
\end{proof}

We shall illustrate below this proposition with the example of Gibbs random fields at
high-temperature under the Dobrushin uniqueness condition.

\section{Examples}

\subsection{High-temperature Gibbs measures}

For the sake of convenience, we briefly recall a few facts about Gibbs measures. We refer to
\cite{Geo} for details.

A finite-range potential (with range $R$) is a family of functions $U(A,\si)$ indexed by finite
subsets $A$ of $\Zd$ such that the value of $U(A,\si)$ depends only on $\si_A$
and such that $U(A,\si)=0$ if $\textup{diam}(A)>R$. If $R=1$ then the potential is nearest-neighbor.

The associated finite-volume Hamiltonian with boundary condition $\eta$ is then given by
$$
H_{\Lambda}^\eta(\si)= \sum_{A\cap \Lambda \neq\emptyset} U(A,\si_{\Lambda}\eta_{\Lambda^c})\,.
$$
The specification is then defined as
$$
\gamma_{\Lambda} (\si|\eta)
=\frac{e^{-H_{\Lambda}^\eta(\si)}}{Z_{\Lambda}^{\eta}}\,\cdot
$$
We then say that $\pee$ is Gibbs measure with potential $U$ if 
$\gamma_{\Lambda} (\si|\cdot)$ is a version of the conditional probability
$\pee(\si_{\Lambda}| \mathcal{F}_{\Lambda^c})$. 

Before we state our result, we need some notions from \cite{ghm}. What we
mean by ``high temperature'' will be an estimate on the variation of
single-site conditional probabilities, which will imply a uniform estimate
for disagreement percolation. For $y\in\Z^d$, let
$$
p_y := 2\sup_{\si,\si'} \left|\pee(\si_y=+|\si_{\Zd\setminus y})
-\pee(\si'_y=+|\si'_{\Zd\setminus y}) \right|\,. 
$$
Writing ${\mathbf p}$ for $(p_y)_y$, let $\nu_{{\mathbf p}}$ denote the
Bernoulli measure on $\{-,+\}^{\Z^d}$ with $\nu_{{\mathbf
p}}(\{X(y)=+\})=p_y$, and $\nu_{p_y}$ its single-site marginal.

From \cite[Theorem 7.1]{ghm} it follows that there exists a coupling
$\pee^{\si}_{x,+,-}$ of the conditional distributions 
$\pee(\cdot | \si_{<x},+_x)$ and $\pee(\cdot | \si_{<x},-_x)$ such that under this coupling
\begin{enumerate}
\item For $y>x$, the event $X_1(y)\neq X_2(y)$ coincides with the event that there exists a
path $\gamma\subset \Zd\setminus(<x)$ from $x$ to $y$ such that, for all $z\in\gamma$, $X_1(z)\neq X_2(z)$.
We denote this event by ``$x\disagree y$''.
\item The distribution of $\1\{X_1(y)\neq X_2(y)\}$ for $y\in\Zd\setminus (\leq x)$ under $\pee^{\si}_{x,+,-}$ 
is dominated by the product measure
$$
\prod_{y\in\Zd\setminus (\leq x)} \nu_{p_y}\,.
$$
\end{enumerate}
Let $p_c=p_c(d)$ be the critical percolation threshold for site-percolation on $\Zd$.
It then follows from statements 1 \& 2 above that, if
\begin{equation}\label{percol}
\sup\{p_y: y\in\Zd\}<p_c
\end{equation}
then we have the uniform estimate
\begin{equation}\label{marilyn}
\pee^{\si}_{x,+,-}\left\{X_1(y)\neq X_2(y)\right\}\leq 
\prod_{y\in\Zd\setminus (\leq x)} \nu_{p_y} (x\disagree y) \leq e^{-c |x-y|}\,.
\end{equation}

Then we can apply Theorem \ref{chardon} to obtain

\begin{theorem}\label{creme}
Let $U$ be a nearest-neighbor potential such that \eqref{percol} holds.
Then for the coupling matrix \eqref{cigare} we have 
the uniform estimate
$$
D_{x,y}^\si \leq e^{-C |x-y|}
$$
for some $C>0$. Hence we have the following exponential concentration inequality:
for any local function $g$ and for all $t>0$
$$
\pee\left\{|g - \E g| \geq t \right\}
\leq 2 \exp\left(-\frac{2 t^2}{\frac{1}{1-e^{-2C}}\  \|\delta g\|_{\ldeuxzd}^2}\right)\cdot
$$
\end{theorem}

\begin{remark}
Theorem \ref{creme} can easily be extended to any finite-range potential. 
\end{remark}

Theorem \ref{creme} was obtained in \cite{Kul} in the Dobrushin's
uniqueness regime \cite[Chapter 8]{Geo}
using a different approach.
The high-temperature condition which we use here is sometimes less
restrictive than Dobrushin's uniqueness condition, 
but sometimes it is more restrictive. However, Dobrushin's uniqueness
condition is not limited to finite-range 
potentials. We now apply Proposition \ref{sardine} to show that in the
Dobrushin's uniqueness regime, there does 
exist a coupling of
$\pee(\cdot | \si_{<x,+_x})$ and $\pee(\cdot | \si_{<x,-_x})$
such that the elements of the associated coupling matrix decay at least as fast as the elements of the
Dobrushin's matrix. 
The Dobrushin's uniqueness condition is based on the matrix
$$
C_{x,y}:= 2 \sup_{\si,\si':\si_{\Zd\setminus y}=\si'_{\Zd\setminus y}} 
\left|\pee(\si_x=+| \si_{\Zd\setminus x})-\pee(\si_x=+ | \si'_{\Zd\setminus x})\right|\,.
$$
This condition is defined by requiring that
$$
\sup_{x\in\Zd} \sum_{y\in\Zd} C_{x,y} <1
$$
and the Dobrushin matrix is then defined as
$$
\Delta_{x,y}:= \sum_{n\geq 0} C_{x,y}^n\,.
$$

We now have the following proposition:

\begin{proposition}
Assume that the Dobrushin uniqueness condition holds.
For any $\varphi:\Zd\to\R^+$ such that for any $x\in\Zd$, 
$$
\sum_{y\in\Zd} \varphi(y) \Delta_{y,x} <\infty\ .
$$
Then there exists a coupling $\hat{\pee}^{\si}_{x,+,-}$ of
$\pee(\cdot | \si_{<x,+_x})$ and $\pee(\cdot | \si_{<x,-_x})$ such that
$$
\sum_{y\in\Zd} \varphi(y)\ \hat{\pee}^{\si_{<x}}_{x,+,-}\left\{ X_1(y)\neq X_2(y)\right\}<\infty\,.
$$
\end{proposition}

\begin{proof}
From \cite[Lemma 1]{Kul}, we have the estimate
$$
\left|
\int g(\eta)\  \pee(d\eta | \si_{<x,+_x})  - \int g(\eta)\ \pee(d\eta | \si_{<x,-_x}) 
\right|
\leq \sum_{y\in\Zd} (\1_{x,y} + \Delta_{y,x}) \delta_y g
$$
(where $\1_{x,y}$ denotes the Kronecker symbol).

We can apply Proposition \ref{sardine} to conclude the proof.
\end{proof}

As an example we mention that if the potential is finite-range and translation-invariant
and satisfies the Dobrushin uniqueness condition, we have
for large enough $|x-y|$
$$
\Delta_{y,x} \leq e^{-c|x-y|}
$$
and hence there exists a coupling $\hat{\pee}^{\si_{<x}}_{x,+,-}$ such that
$$
\hat{\pee}^{\si_{<x}}_{x,+,-}\left\{ X_1(y)\neq X_2(y)\right\} \leq e^{-c'|x-y|}
$$
for all $c'<c$ and large enough $|x-y|$.

Unfortunately, we are not able to construct explicitly such a coupling.

\subsection{The low-temperature Ising model}

It is clear that for the Ising model in the phase
coexistence region, no exponential concentration inequalities
can hold. Indeed, this would contradict the surface-order large deviation bounds for
the magnetization in that regime (see e.g. \cite{ioffe} and references therein).
Nevertheless, we shall show that we can control all moments and obtain
stretched-exponential inequalities (which are compatible with large
deviation bounds).

We consider the low-temperature plus phase
of the Ising model on $\Z^d$, $d\geq 2$. This is a probability measure
$\pee^+_\beta$ on lattice spin configurations
$\si\in\Omega $, defined as the weak limit as
$\Lambda\uparrow\Z^d$
of the following finite-volume measures:
\begin{equation}\label{rat}
\pee^+_{\Lambda,\beta}(\si_\Lambda) =
\exp\left(\beta\sum_{<xy>\in \Lambda} \si_x\si_y + \beta\sum_{<xy>, x\in\partial
\Lambda, \ y\notin \Lambda}\si_x\right)
\Big/ Z_{\Lambda,\beta}^+
\end{equation}
where $\beta\in\R^+$ is the inverse temperature, and
$Z_{\Lambda,\beta}^+$ is the partition function. 
In \eqref{rat} $<xy>$ denotes nearest neighbor bonds and $\partial \Lambda$ the
inner boundary,
i.e. the set of those $x\in \Lambda$ having at least one neighbor
$y\notin \Lambda$. 
The existence of the limit $\Lambda\uparrow \Z^d$ of
$\pee_{\Lambda,\beta}^+$ is by a standard and well-known 
monotonicity argument, see
e.g. \cite{Geo}.

For any $\eta\in \Omega$ $, \Lambda\subset \Z^d$ we denote by
$\pee_{\Lambda,\beta}^\eta$ the corresponding finite-volume
measure with boundary condition $\eta$:
$$
\pee^\eta_{\Lambda,\beta}(\si_\Lambda) = \exp\left(\beta \sum_{<xy>\in \Lambda}
\si_x\si_y+ \beta \sum_{ x\in \Lambda,\ y\notin
\Lambda}\si_x\eta_x\right)
\Big/ Z_{\Lambda,\beta}^{\eta}\,.
$$
Later on we will have to choose $\beta$ large enough, in particular,
greater than the critical inverse temperature $\beta_c$ ($\beta<\beta_c$
implies uniqueness of the infinite-volume measure).

We can now formulate our results on arbitrary local functions for the
low-temperature Ising model.

\begin{theorem}
Let $\pee=\pee_\beta^+$ be the plus phase of the low-temperature Ising
model defined above. There exists $\beta_0>\beta_c$, such that for all
$\beta>\beta_0$, for any local function $g$, we have the following
inequalities:
\begin{itemize}
\item For all $p\in\N$, there exists a constant $C_p\in ]0,\infty[$ such that 
$$
\E\left[(g-\E g)^{2p}\right] \leq C_p\ \|\delta g\|_{\ldeuxzd}^{2p}\,.
$$
Consequently, for all $t>0$, we have the concentration inequalities
$$
\pee\left\{ \left|g-\E g \right|>t\right\} \leq 
C_p\ \frac{\| \delta g\|_{\ldeuxzd}^{2p}}{t^{2p}}\cdot
$$
\item 
Moreover, there exists $0<\varrho(\beta)<1$, such that for any 
$0<\varrho<\varrho(\beta)$ there is a constant $K_{\varrho}>0$, such that 
we have, for any local function $g$,
$$
\left\|g-\E g\right\|_{\Phi_{\varrho}}
\le K_{\varrho}\ \|\delta g\|_{\ldeuxzd}\,.
$$
Consequently, there exists a constant $c_\varrho\in]0,\infty[$
such that, for all $t>0$,
$$
\pee\left\{ |g -\E g| >t\right\} \leq 4
\exp\left(-c_\varrho\ \frac{t^\varrho}{\|\delta g\|_{\ldeuxzd}^\varrho}\right)\cdot
$$
\end{itemize}
\end{theorem}

\begin{proof}
This theorem is an application of Theorem \ref{curantobis} and
Proposition \ref{pinardbis}. All we have to do is to 
obtain the bound \eqref{lizd} with good decay properties for
the tail of the distribution of $\ell_0$ to ensure the finiteness
of $\E(\ell_0^{2pd+\epsilon})$, $\|\ell_0^d\|_{\Phi_\vartheta}$, and
of $\|\psi\|_{\lun}$. This is an immediate consequence of the next proposition.
\end{proof}

\begin{proposition}
Let $\pee=\pee_{\beta}^+$ be the plus phase of the low-temperature
Ising model.
There exists $\beta_0>\beta_c$ such that for all $\beta>\beta_0$, the inequality \eqref{lizd} holds
together with the estimate
$$
\psi(n)\leq C e^{-c n}
$$
for all $n\in\N$ and
$$
\pee\{\ell_0 \geq n\} \leq C' e^{-c' n^{\alpha}}
$$
for some $c,c',C,C'>0$ and $0<\alpha\leq 1$.
\end{proposition}

\begin{proof}
We shall make a coupling of the conditional measures
$\pee(\cdot | \si_{<x,+_x})$ and $\pee(\cdot | \si_{<x,-_x})$.
This coupling already appeared in \cite{vdbm} (see also \cite{ghm}).
Both conditional measures are a distribution of a random field $\omega_y$, $y\notin (\leq x)$.
We start with the first site $y_1>x$ according to the order induced by
$\Gamma$ (see Section \ref{RF}).
We generate $X_1(y_1)$ and $X_2(y_1)$ as a realization of the maximal coupling between
$\pee(\sigma_{y_1}=\cdot | \si_{<x,+_x})$ and $\pee(\sigma_{y_1}=\cdot | \si_{<x,-_x})$.
Given that we have generated $X_1(y),X_2(y),\ldots,X_1(y_n),X_2(y_n)$ for $y=y_1,\ldots,y_n$, we generate
$X_1(y_{n+1})$, $X_2(y_{n+1})$ for the smallest $y_{n+1}>y_n$ as a realization of the maximal
coupling between
$$
\pee(\sigma_{y_{n+1}}=\cdot | X_1(y_1)\cdots X_1(y_n)\si_{<x,+_x})\;\textup{and}\;
\pee(\sigma_{y_{n+1}}=\cdot | X_2(y_1)\cdots X_2(y_n)\si_{<x,-_x})\,.
$$
By the Markov property of $\pee$ we have the following:
if there exists a contour separating $y$ from $x$ such that for all sites $z$ belonging to
that contour we have $X_1(z)=X_2(z)$, then $X_1(y)=X_2(y)$.
The complement of this event (of having such a contour) is contained in the event that there
exists a path of disagreement from $x$ to $y$, i.e., a path $\gamma\subset \Zd\setminus (<x)$ 
such that for all $z\in\gamma$, $X_1(z)\neq X_2(z)$. Denote that event by $E_{xy}$. Clearly its
probability is bounded from above by the probability of the same event in the product coupling.
In turn the event $E_{xy}$ is contained in the event $E^+_{xy}$ that there exists a path $\gamma$
from $x$ to $y$ in $\Zd\setminus (<x)$ such that for all $z\in\gamma$, $(X_1(z),X_2(z))\neq (+,+)$.
In \cite{mrsv} the probability of that event in the product coupling is precisely estimated from
above by
\begin{equation}\label{himmler}
C e ^{-c |x-y|} + \1\{ \ell_x (\si) \geq |x-y|\}
\end{equation}
for some $C,c>0$, where $\ell_x (\si)$ is an unbounded function of $\si$ with tail estimate
$$
\pee(\ell_x (\si)\geq n)=\pee(\ell_0 (\si)\geq n)
\leq C' e^{-c'n^{\alpha}} 
$$
for some $C',c'>0$ and $0<\alpha<1$.
For the reader's convenience, we briefly comment on these estimates. The idea is that the conditional
measure $\pee(\cdot | \xi_{\leq x})$ resembles the original unconditioned plus phase (in
$\Zd\setminus (\leq x)$) provided $\xi$ contains ``enough'' pluses. ``Containing enough pluses''
is exactly quantified by the random variable $\ell_x(\xi)$: $(\ell_x(\xi)\leq n)$
is the event that for all self-avoiding path $\gamma$ of length at least $n$, the magnetization
along $\gamma$, 
$$
m_\gamma(\xi):= \frac{1}{|\gamma|} \sum_{z\in \gamma} \xi_z
$$
is close ``enough to one''. If this is the case then under the conditional measure we still
have a Peierls' estimate, which produces the exponential term in \eqref{himmler}.
We refer to \cite{mrsv} for more details.
\end{proof}

\bigskip

\noindent {\bf Acknowledgement}. We wish to thank a referee for a
careful reading of our manuscript, as well as for the impetus he gave
us to include  stretched-exponential concentration inequalities which we did
not consider in the primitive version of this work.


\end{document}